\documentclass[a4paper]{article}
\usepackage{amsmath,amssymb,mathrsfs,graphicx,psfrag}
\DeclareMathOperator{\re}{Re}
\DeclareMathOperator{\im}{Im}
\DeclareMathOperator{\fp}{f.p.}
\newtheorem{thm}{Theorem}
\allowdisplaybreaks[1]
\begin{document}
\title{A numerical method of computing Hadamard finite-part integrals with an integral power singularity 
at the endpoint on a half infinite interval}
\author{Hidenori Ogata%
\footnote{Department of Computer and Network Engineering, %
Graduate School of Informatics and Engineering, %
The University of Electro-Communications, 5
1-5-1 Chofugaoka, Chofu, Tokyo 182-8585, Japan, %
(e-mail) {\tt ogata@im.uec.ac.jp}}}
\maketitle
\begin{abstract}
 In this paper, we propose a numerical method of computing Hadamard finite-part integrals with an integral 
 power singularity at the endpoint on a half infinite interval, that is, a finite value assigned to 
 a divergent integral with an integral power singularity at the endpoint on a half infinite interval. 
 In the proposed method, we express a desired finite-part integral using a complex integral, and we obtain 
 the integral by evaluating the complex integral by the DE formula. 
 Theoretical error estimate and some numerical examples show the effectiveness of the proposed method. 
\end{abstract}
%%%%%%%%%%%%%%%%%%%%%%%%%%%%%%%%%%%%%%%%%%%%%%%%%%%%%%%%%%%%%%%%%%%%%%%%%%%%%%%%%%%%%%%%%%%%%%%%%%%%%
\section{Introduction}
\label{sec:introduction}
The integral 
\begin{equation*}
 \int_0^{\infty}x^{-1}f(x)\mathrm{d}x, 
\end{equation*}
where $f(x)$ is an analytic function on the half infinite interval $[0,+\infty)$ such that 
$f(0)\neq 0$ and $f(x) = \mathrm{O}(x^{-\alpha})$ with $\alpha>0$, is divergent. 
However, we can assign a finite value to this divergent integral. 
In fact, we consider the integral
\begin{equation*}
 \int_{\epsilon}^{\infty}x^{-1}f(x)\mathrm{d}x
\end{equation*}
with $\epsilon > 0$, and, by integrating by part, we have
\begin{align*}
 \int_{\epsilon}^{\infty}x^{-1}f(x)\mathrm{d}x = \: & 
 \int_{\epsilon}^{\infty}(\log x)^{\prime}f(x)\mathrm{d}x
 \\ 
 = \: & 
 \bigg[\log x f(x)\bigg]_{\epsilon}^{\infty} - 
 \int_{\epsilon}^{\infty}\log x f^{\prime}(x)\mathrm{d}x
 \\ 
 = \: & 
 - f(\epsilon)\log\epsilon - \int_{\epsilon}^{\infty}\log x f^{\prime}(x)\mathrm{d}x
 \\ 
 = \: & 
 - f(0)\log\epsilon - \int_{0}^{\infty}\log x f^{\prime}(x)\mathrm{d}x 
 + \mathrm{O}(1) \quad ( \: \mbox{as} \ \epsilon\downarrow 0 \: ).
\end{align*}
Then, the limit
\begin{equation*}
 \lim_{\epsilon\downarrow 0}
  \left\{
   \int_{\epsilon}^{\infty}x^{-1}f(x)\mathrm{d}x + f(0)\log\epsilon
  \right\}
\end{equation*}
exists and is finite. 
We call this limit the Hadamard finite-part (f.p.) integral and denote it by 
\begin{equation*}
 \fp\int_0^{\infty}x^{-1}f(x)\mathrm{d}x.
\end{equation*}
In general, we can define the f.p. integral
\begin{equation}
 \label{eq:fp-integral0}
  I^{(n)}[f] = \fp\int_0^{\infty}x^{-n}f(x)\mathrm{d}x \quad ( \: n = 1, 2, \ldots \: ), 
\end{equation}
where $f(x)$ is an analytic function such that $f(0)\neq 0$  
and $f(x)=\mathrm{O}(x^{n-\alpha-1})$ $(\:x\rightarrow +\infty\:)$ with $\alpha>0$ 
\cite{EstradaKanwal1989}. 

In this paper, we propose a numerical method of computing f.p. integrals (\ref{eq:fp-integral0}). 
In the proposed method, we express the f.p. integral (\ref{eq:fp-integral0}) using a complex integral, 
and we obtain the f.p. integral by evaluating the complex integral by the DE formula \cite{TakahasiMori1974}. 
Theoretical error estimate and some numerical examples show the exponential convergence of the proposed formula 
as the number of sampling points increases. 

Previous works related to this paper are as follows. 
The author and Hirayama proposed a numerical method of computing ordinary integrals related to hyperfunction theory 
\cite{OgataHirayama2018}, a theory of generalized functions based on complex function theory. 
The author proposed numerical methods for computing Hadamard finite-part integrals with a singularity at an endpoint 
on a finite interval \cite{Ogata2019b,Ogata2019c}. 
In these methods, we express a desired integral using a complex integral, we obtain the integral by evaluating the 
complex integral by conventional numerical integration formulas. 
For Cauchy principal-value integrals or Hadamard finite-part integrals on a finite interval with a singularity 
in the interior of the integral interval
\begin{equation}
 \label{eq:fp-integral02}
 \fp\int_0^1 \frac{f(x)}{(x-\lambda)^n}\mathrm{d}x \quad 
  ( \: 0 < \lambda < 1, \: n = 1, 2, \ldots \: ),
\end{equation}
many methods were proposed. 
Elliot and Paget proposed Gauss-type numerical integration formulas for (\ref{eq:fp-integral02})  
\cite{ElliotPaget1979,Paget1981}. 
Bialecki proposed Sinc numerical integration formulas for (\ref{eq:fp-integral02}) 
\cite{Bialecki1990a,Bialecki1990b}, 
where the trapezoidal formula is used together with variable transform technique as in the DE formula 
\cite{TakahasiMori1974}. 
Ogata and et al. improved them and proposed a DE-type numerical integration formula for 
(\ref{eq:fp-integral02}) \cite{OgataSugiharaMori2000}.

The remainder of this paper is structured as follows. 
In Section \ref{sec:fp-integral}, we define the f.p. integral (\ref{eq:fp-integral0}) and propose 
a numerical method of computing it. 
In addition, we show theoretical error estimate of the proposed method. 
In Section \ref{sec:example}, we show some numerical example which show the effectiveness of the proposed method. 
In Section \ref{sec:summary}, we give a summary of this paper. 
%%%%%%%%%%%%%%%%%%%%%%%%%%%%%%%%%%%%%%%%%%%%%%%%%%%%%%%%%%%%%%%%%%%%%%%%%%%%%%%%%%%%%%%%%%%%%%%%%%%%%
\section{Hadamard finite-part integrals and a numerical method}
\label{sec:fp-integral}
The f.p. integral (\ref{eq:fp-integral0}) is defined by
\begin{multline}
 \label{eq:fp-integral}
 I^{(n)}[f] = \fp\int_0^{\infty}x^{-n}f(x)\mathrm{d}x 
 \\ 
 = 
 \lim_{\epsilon\downarrow 0}\left\{
 \int_{\epsilon}^{\infty}x^{-n}f(x)\mathrm{d}x 
 - 
 \sum_{k=0}^{n-2}\frac{\epsilon^{k+1-n}}{k!(n-1-k)}f^{(k)}(0) 
 + 
 \frac{\log\epsilon}{(n-1)!}f^{(n-1)}(0)
 \right\}
 \\
 ( \: n = 1, 2, \ldots \: ),
\end{multline}
where $f(x)$ is an analytic function on $[0,+\infty)$ such that $f(0)\neq 0$ and 
$f(x) = \mathrm{O}(x^{n-1-\alpha})$ as $x\rightarrow +\infty$ with $\alpha>0$, 
and the second term on the right-hand side is zero if $n=1$. 
We can show that it is well-defined as follows. 
In fact, for $\epsilon>0$, we can show by integrating by part
\begin{align*}
 & \int_{\epsilon}^{\infty}x^{-n}f(x)\mathrm{d}x
 \\ 
 = \: & 
 \frac{\epsilon^{1-n}}{n-1}f(\epsilon) 
 + 
 \frac{1}{n-1}\int_{\epsilon}^{\infty}x^{1-n}f^{\prime}(x)\mathrm{d}x
 \\ 
 = \: & 
 \frac{\epsilon^{1-n}}{n-1}f(\epsilon) + 
 \frac{\epsilon^{2-n}}{(n-1)(n-2)}f^{\prime}(\epsilon) 
 + 
 \frac{1}{(n-1)(n-2)}\int_{\epsilon}^{\infty}x^{2-n}f^{\prime\prime}(x)\mathrm{d}x
 \\ 
 = \: & \cdots
 \\ 
 = \: & 
 \sum_{k=0}^{n-2}\frac{\epsilon^{k+1-n}}{(n-1)(n-2)\cdots(n-1-k)}f^{(k)}(\epsilon)
 - 
 \frac{\log\epsilon}{(n-1)!}f^{(n-1)}(\epsilon)
 \\ 
 & + 
 \frac{1}{(n-1)!}\int_{\epsilon}^{\infty}\log x f^{(n)}(x)\mathrm{d}x
 \\ 
 = \: & 
 \sum_{k=0}^{n-2}\frac{\epsilon^{k+1-n}}{(n-1)(n-2)\cdots(n-1-k)}
 \left\{\sum_{l=0}^{\infty}\frac{f^{(k+l)}(0)}{l!}\epsilon^l\right\}
 \\ 
 & 
 - 
 \frac{\log\epsilon}{(n-1)!}\sum_{k=0}^{\infty}\frac{f^{(k)}(0)}{k!}\epsilon^k
 + 
 \frac{1}{(n-1)!}\int_{\epsilon}^{\infty}\log x f^{(n)}(x)\mathrm{d}x
 \\ 
 = \: & 
 \sum_{l=0}^{n-2}
 \bigg\{
 \underbrace{%
 \sum_{k=0}^{l}\frac{1}{(l-k)!(n-1)(n-2)\cdots(n-1-k)}%
 }_{(\ast)}
 \bigg\}\epsilon^{l+1-n}f^{(l)}(0)
 \\ 
 & 
 - \frac{\log\epsilon}{(n-1)!}f^{(n-1)}(0) + \mathrm{O}(1) 
 \quad ( \: \mbox{as} \ \epsilon\downarrow 0 \: ),
\end{align*}
and
\begin{align*}
 (\ast) = \: & 
 \frac{1}{l!(n-1)} + \frac{1}{(l-1)!(n-1)(n-2)}
 \\ 
 & 
 + \cdots 
 + \frac{1}{1!(n-1)(n-2)\cdots(n-l+1)(n-l)} 
 \\ 
 & 
 + \frac{1}{(n-1)(n-2)\cdots(n-l+1)(n-l)(n-l-1)}
 \\
 = \: & 
 \frac{1}{l!(n-1)} + \frac{1}{(l-1)!(n-1)(n-2)} + \cdots 
 \\ 
 & 
 + \frac{1}{2!(n-1)(n-2)\cdots(n-l+2)(n-l+1)} 
 \\ 
 & 
 + \frac{1}{1!(n-1)(n-2)\cdots(n-l+2)(n-l+1)(n-l-1)}
 \\
 = \: & 
 \frac{1}{l!(n-1)} + \frac{1}{(l-1)!(n-1)(n-2)} + \cdots 
 \\ 
 & 
 + \frac{1}{3!(n-1)(n-2)\cdots(n-l+3)(n-l+2)} 
 \\ 
 & 
 + \frac{1}{2!(n-1)(n-2)\cdots(n-l+3)(n-l+2)(n-l-1)}
 \\ 
 = \: & \cdots 
 \\ 
 = \: & 
 \frac{1}{l!(n-l-1)}.
\end{align*}
Then, we have
\begin{align*}
 \int_{\epsilon}^{\infty}x^{-n}f(x)\mathrm{d}x = \: & 
 \sum_{l=0}^{n-2}\frac{\epsilon^{l+1-n}}{l!(n-l-1)}f^{(l)}(0) 
 - 
 \frac{\log\epsilon}{(n-1)!}f^{(n-1)}(0) 
 \\ 
 & + \mathrm{O}(1)
 \quad ( \: \mbox{as} \ \epsilon\downarrow 0 \: ).
\end{align*}
Therefore, the limit in (\ref{eq:fp-integral}) exists and is finite.

The f.p. integral (\ref{eq:fp-integral}) is expressed using a complex integral.
\begin{thm}
 \label{thm:complex-integral}
 We suppose that $f(z)$ is analytic in a complex domain $D$, which contains the half infinite interval 
 $[0,+\infty)$ in its interior. Then, the f.p. integral (\ref{eq:fp-integral}) is expressed as 
 \begin{equation}
  \label{eq:complex-integral}
   I^{(n)}[f] = 
   %   \fp\int_0^{\infty} x^{-n}f(x)\mathrm{d}x = 
   \frac{1}{2\pi\mathrm{i}}\oint_C z^{-n}f(z)\log(-z)\mathrm{d}z, 
 \end{equation}
 where $C$ is a complex integral path such that it encircles $[0,+\infty)$ in the positive sense and 
 is contained in $D$. 
\end{thm}

\begin{figure}[htbp]
 \begin{center}
  \psfrag{0}{$\mathrm{O}$}
  \psfrag{a}{$C_{\epsilon}$}
  \psfrag{e}{$\epsilon$}
  \psfrag{p}{$\Gamma_{\epsilon}^{(+)}$}
  \psfrag{m}{$\Gamma_{\epsilon}^{(-)}$}
  \includegraphics[width=0.7\textwidth]{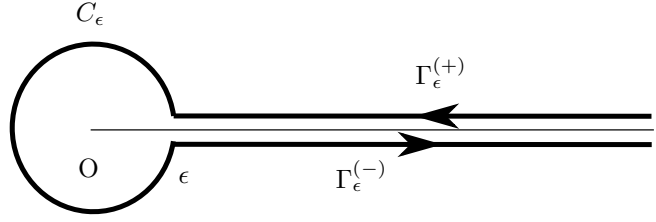}
 \end{center}
 \caption{The complex integral paths $\Gamma_{\epsilon}^{(\pm)}$ and $C_{\epsilon}$.}
 \label{fig:proof-integral-path}
\end{figure}

\paragraph{Proof of Theorem \ref{thm:complex-integral}}
From Cauchy's integral theorem, we have
\begin{equation*}
 \frac{1}{2\pi\mathrm{i}}\oint_C z^{-n}f(z)\log(-z)\mathrm{d}z
  = 
  \frac{1}{2\pi\mathrm{i}}
  \left(\int_{\Gamma_{\epsilon}^{(+)}}+\int_{C_{\epsilon}}+\int_{\Gamma_{\epsilon}^{(-)}}\right)
  z^{-n}f(z)\log(-z)\mathrm{d}z, 
\end{equation*}
for $\epsilon>0$, where $\Gamma_{\epsilon}^{(+)}$ and $C_{\epsilon}$ are complex integral paths respectively 
defined by
\begin{align*}
 \Gamma_{\epsilon}^{(+)} = \: & 
 \{ \: x\in\mathbb{R}+\mathrm{i}0 \: | \: +\infty > x \geqq \epsilon \: \}, 
 \\
 \Gamma_{\epsilon}^{(-)} = \: & 
 \{ \: x\in\mathbb{R}-\mathrm{i}0 \: | \: \epsilon \leqq x < +\infty \: \}, 
 \\ 
 C_{\epsilon} = \: & 
 \{ \: \epsilon\mathrm{e}^{\mathrm{i}\theta} \: | \: 0 \leqq \theta \leqq 2\pi \: \}
\end{align*}
(see Figure \ref{fig:proof-integral-path}), 
and the complex logarithmic function $\log z$ is the principal value, that is, the branch such that 
it takes a real value on the positive real axis. 
As to the integrals on $\Gamma_{\epsilon}^{(\pm)}$, we have
\begin{align*}
 & 
 \frac{1}{2\pi\mathrm{i}}
 \left( \int_{\Gamma_{\epsilon}^{(+)}} + \int_{\Gamma_{\epsilon}^{(-)}} \right)
 z^{-n}f(z)\log(-z)\mathrm{d}z
 \\
 = \: & 
 \frac{1}{2\pi\mathrm{i}}
 \left\{
 - \int_{\epsilon}^{\infty}x^{-n}f(x)\log(x\mathrm{e}^{-\mathrm{i}\pi})\mathrm{d}x
 + \int_{\epsilon}^{\infty}x^{-n}f(x)\log(x\mathrm{e}^{\mathrm{i}\pi})\mathrm{d}x
 \right\}
 \\ 
 = \: & 
 \frac{1}{2\pi\mathrm{i}}
 \left\{
 - \int_{\epsilon}^{\infty}x^{-n}f(x)(\log x - \mathrm{i}\pi)\mathrm{d}x
 + \int_{\epsilon}^{\infty}x^{-n}f(x)(\log x + \mathrm{i}\pi)\mathrm{d}x
 \right\}
 \\ 
 = \: & 
 \int_{\epsilon}^{\infty}x^{-n}f(x)\mathrm{d}x.
\end{align*}
As to the integral on $C_{\epsilon}$, we have
\begin{align}
 \nonumber
 & 
 \frac{1}{2\pi\mathrm{i}}\int_{C_{\epsilon}}z^{-n}f(z)\log(-z)\mathrm{d}z
 \\ 
 \nonumber
 = \: & 
 \frac{1}{2\pi\mathrm{i}}\int_0^{2\pi}(\epsilon\mathrm{e}^{\mathrm{i}\theta})^{-n}
 f(\epsilon\mathrm{e}^{\mathrm{i}\theta})\log(\epsilon\mathrm{e}^{\mathrm{i}(\theta-\pi)})
 \mathrm{i}\epsilon\mathrm{e}^{\mathrm{i}\theta}\mathrm{d}\theta
 \\ 
 \nonumber
 = \: & 
 \frac{\epsilon^{1-n}}{2\pi}\int_{0}^{2\pi}\mathrm{e}^{\mathrm{i}(1-n)\theta}
 \left\{\sum_{k=0}^{\infty}\frac{f^{(k)}(0)}{k!}\epsilon^k\mathrm{e}^{\mathrm{i}k\theta}\right\}
 \left\{\log\epsilon + \mathrm{i}(\theta-\pi)\right\}
 \mathrm{d}\theta
 \\
 \nonumber
 = \: & 
 \frac{1}{2\pi}\sum_{k=0}^{\infty}\frac{\epsilon^{k-n+1}}{k!}\log\epsilon f^{(k)}(0)
 \int_0^{2\pi}\mathrm{e}^{\mathrm{i}(k-n+1)\theta}\mathrm{d}\theta
 \\ 
 \label{eq:proof-complex-integral}
 & 
 + 
 \frac{\mathrm{i}}{2\pi}\sum_{k=0}^{\infty}\frac{\epsilon^{k-n+1}}{k!}f^{(k)}(0)
 \int_0^{2\pi}(\theta-\pi)\mathrm{e}^{\mathrm{i}(k-n+1)\theta}
 \mathrm{d}\theta,
\end{align}
where we exchanged the order of the integral and the infinite sum since the infinite series is 
uniformly convergent on $0\leqq\theta\leqq 2\pi$.
Since
\begin{gather*}
 \int_{0}^{2\pi}\mathrm{e}^{\mathrm{i}(k-n+1)\theta}\mathrm{d}\theta = 
 2\pi\delta_{k,n-1}, 
 \\
 \int_{0}^{2\pi}(\theta-\pi)\mathrm{e}^{\mathrm{i}(k-n+1)\theta}\mathrm{d}\theta 
 = 
 \begin{cases}
  \displaystyle \frac{-2\pi}{\mathrm{i}(n-1-k)} & ( \: 0 \leqq k \leqq n-2 \: ) \\ 
  0 & ( \: k = n-1 \: ),
 \end{cases}
\end{gather*}
we have
\begin{equation*}
 \mbox{(\ref{eq:proof-complex-integral})}  
 = 
 \frac{\log\epsilon}{(n-1)!}f^{(n-1)}(0) - 
 \sum_{k=0}^{n-2}\frac{\epsilon^{k-n+1}}{k!(n-1-k)}f^{(k)}(0) + 
 \mathrm{O}(\epsilon) \quad ( \: \epsilon\downarrow 0 \: ).
\end{equation*}
Summarizing the above calculations, we have
\begin{align*}
 \frac{1}{2\pi\mathrm{i}}\oint_C z^{-n}f(z)\log(-z)\mathrm{d}z
 = \: & 
 \int_{\epsilon}^{1}x^{-n}f(x)\mathrm{d}x 
 - 
 \sum_{k=0}^{n-2}\frac{\epsilon^{k-n+1}}{k!(n-1-k)}f^{(k)}(0)
 \\
 & 
 + 
 \frac{\log\epsilon}{(n-1)!}f^{(n-1)}(0) 
 + \mathrm{O}(\epsilon) \quad ( \: \epsilon\downarrow 0 \: ), 
\end{align*}
and, taking the limit $\epsilon\downarrow 0$, we obtain (\ref{eq:complex-integral}).
\hfill\rule{1.5ex}{1.5ex}

We obtain the f.p. integral by evaluating the complex integral on the right-hand side of (\ref{eq:complex-integral}) 
by a conventional numerical integration formula such as the DE formula \cite{TakahasiMori1974}, that is, 
\begin{equation*}
 \int_{-\infty}^{\infty}g(u)\mathrm{d}u = 
 \int_{-\infty}^{\infty}g(\psi_{\rm DE}(v))\psi_{\rm DE}^{\prime}(v)\mathrm{d}v
 \simeq 
 h\sum_{k=-N_-}^{N_+}g(\psi_{\rm DE}(kh))\psi_{\rm DE}^{\prime}(kh), 
\end{equation*}
where $h>0$ is the mesh of the trapezoidal formula, 
$u = \psi_{\rm DE}(v)$ is the DE transform
\begin{equation*}
 \psi_{\rm DE}(v) = 
  \begin{cases}
   \sinh(\sinh v) & 
   ( \: g(u) = \mathrm{O}(|u|^{-\alpha-1}) \quad \mbox{as} \quad u\rightarrow\pm\infty, \: \alpha > 0 \: )
   \\
   \sinh v & ( \: g(u) = \mathrm{O}(\exp(-c|u|)) \quad \mbox{as} \quad u\rightarrow\pm\infty, \: c > 0 \: ), 
  \end{cases}
\end{equation*}
and $N_{\pm}$ is a positive integer such that the transformed integrand \\
$g(\psi_{\rm DE}(kh))\psi_{\rm DE}^{\prime}(kh)$ is sufficiently small at $k=-N_-, N_+$. 
We can take $N_{\pm}$ small since $g(\psi_{\rm DE}(v))\psi_{\rm DE}^{\prime}(v)$ decays double exponentially 
as $v\rightarrow\pm\infty$. 
Then, we have the approximation formula
\begin{align}
 % \nonumber
 % & 
 % \fp\int_0^{\infty}x^{-n}f(x)\mathrm{d}x
 % \\
 % \nonumber
 % = \: & 
 % \frac{1}{2\pi\mathrm{i}}\int_{-\infty}^{\infty}\varphi(u)^{-n}
 % f(\varphi(u))\log(-\varphi(u))\varphi^{\prime}(u)\mathrm{d}u
 % \\ 
 \nonumber  
 I^{(n)}[f] \simeq \: & I_{h,N_+,N_-}^{(n)}[f]
 \\ 
 \nonumber
 \equiv \: & 
 \frac{h}{2\pi\mathrm{i}}\sum_{k=-N_-}^{N_+}
 \varphi(\psi_{\rm DE}(kh))^{-n}f(\varphi(\psi_{\rm DE}(kh)))
 \log(-\varphi(\psi_{\rm DE}(kh)))
 \\ 
 \label{eq:approx-fp-integral}
 & \hspace{15mm}
 \times \varphi^{\prime}(\psi_{\rm DE}(kh))
 \psi_{\rm DE}^{\prime}(kh),
\end{align}
where $z = \varphi(u), \ -\infty < u < +\infty$ is a parameterization of the complex integral path $C$. 

If $f(x)$ is an analytic function on the real axis and $C$ is an analytic curve, 
the proposed approximation (\ref{eq:approx-fp-integral}) converges exponentially as shown in the following theorem. 
For the simplicity, we take $N_+=N_-\equiv N^{\prime}$. 
\begin{thm}
 \label{thm:error-estimate}
 We suppose that 
 \begin{enumerate}
  \item the parameterization function $\varphi(w)$ of $C$ is analytic in the strip 
	\begin{equation*}
	 \mathscr{D}_d = \{ \: w\in\mathbb{C} \: | \: |\im w| < d \: \} \quad ( \: d > 0 \: ) 
	\end{equation*}
	such that 
	\begin{equation*}
	 \varphi(\mathscr{D}_d) = \{ \: \varphi(w) \: | \: w\in\mathscr{D}_d \: \}, 
	\end{equation*}
	is contained in $\mathbb{C}\setminus [0,+\infty]$, 
  \item 
	\begin{equation*}
	 \begin{aligned}
	  & \mathscr{N}(f,\varphi,\psi_{\rm DE}, \mathscr{D}_d) 
	  \\
	  \equiv \: & 
	  \lim_{\epsilon\rightarrow 0} 
	  \oint_{\partial\mathscr{D}_d(\epsilon)}
	  |\varphi(\psi_{\rm DE}(w))^{-n}f(\varphi(\psi_{\rm DE}(w)))
	  \log(-\varphi(\psi_{\rm DE}(w)))\psi_{\rm DE}^{\prime}(w)| 
	  \\ 
	  < \: & \infty, 
	 \end{aligned}
	\end{equation*}
	where
	\begin{equation*}
	 \mathscr{D}_d(\epsilon) \equiv 
	  \left\{ \: w \in \mathbb{C} \: | \: |\re w| < 1/\epsilon, \quad |\im w| < d(1-\epsilon) \: \right\}. 
	\end{equation*}
	and 
  \item there exist positive numbers $C_0$, $c_1$ and $c_2$ such that
	\begin{equation*}
	 |f(\varphi(\psi_{\rm DE}(v)))| \leqq C_0\exp(-c_1\exp(c_2|v|)) \quad 
	  ( \: \forall v \in \mathbb{R} \: ). 
	\end{equation*}
 \end{enumerate}
 Then, we have the inequality
 \begin{align}
  \nonumber
  |I^{(n)}[f] - I_{h,N^{\prime}}^{(n)}[f]|
  %  \left|\fp\int_0^{\infty}x^{-n}f(x)\mathrm{d}x - I_{h,N^{\prime}}^{(n)}[f]\right| 
  \leqq \: & 
  \frac{1}{2\pi}\mathscr{N}(f, \varphi, \psi_{\rm DE}, \mathscr{D}_d)
  \frac{\exp(-2\pi d/h)}{1-\exp(-2\pi d/h)}
  \\
  \label{eq:error-estimate}
  & 
  + 
  C(f,\varphi,\psi_{\rm DE}, \mathscr{D}_d)\exp(-c_1\exp(c_2N^{\prime}h)),
 \end{align}
 where $I_{h,N^{\prime}}^{(n)}[f] = I_{h,N^{\prime},N^{\prime}}^{(n)}[f]$ and 
 $C(f, \varphi, \psi_{\rm DE}, \mathscr{D}_d)$ is a positive number depending on 
 $f(z)$, $\varphi$, $\psi_{\rm DE}$ and $\mathscr{D}_d$ only. 
\end{thm}
This theorem shows that the approximation formula (\ref{eq:approx-fp-integral}) converges 
exponentially as the mesh $h$ decreases and the number of sampling points $2N^{\prime}+1$ increases. 
\paragraph{Proof of Theorem \ref{thm:error-estimate}}
We have
\begin{align}
 \nonumber
 & 
 \left|\fp\int_0^{\infty}x^{-n}f(x)\mathrm{d}x - I_{h,N}^{(n)}[f]\right|
 \\
 \nonumber
 \leqq \: & 
 \left|\fp\int_0^{\infty}x^{-n}f(x)\mathrm{d}x - I_h^{(n)}[f] \right| 
 \\ 
 \label{eq:proof-error-estimate}
 & 
 + 
 h\sum_{|k|>N}
 \left|
 \varphi(\psi_{\rm DE}(kh))^{-n}f(\varphi(\psi_{\rm DE})(kh))
 \varphi^{\prime}(\psi_{\rm DE}(kh))\psi_{\rm DE}^{\prime}(kh)
 \right|,
\end{align}
where $I_h^{(n)}[f] = \lim_{N\rightarrow\infty}I_{h,N}^{(n)}[f]$. 
For the first term on the right-hand side of (\ref{eq:proof-error-estimate}), we have
\begin{equation*}
 \left|\fp\int_0^{\infty}x^{-n}f(x)\mathrm{d}x - I_h^{(n)}[f]\right| \leqq 
  \frac{1}{2\pi}\mathscr{N}(f,\varphi,\psi_{\rm DE}, \mathscr{D}_d)\frac{\exp(-2\pi d/h)}{1-\exp(-2\pi d/h)}
\end{equation*}
by Theorem 3.2.1 in \cite{Stenger1993}. 
For the second term on the right-side hand, we have
\begin{align*}
 |\mbox{the second term}| \leqq \: & 
 C_0 h\sum_{|k|>N}\exp(-c_1\exp(c_2 kh))
 \\
 \leqq \: & 
 2C_0 \int_{kh}^{\infty}\exp(-c_1\exp(c_2 x))\mathrm{d}x
 \\
 \leqq \: & 
 2C_0 \int_{kh}^{\infty}\exp(c_2 x)\exp(-c_1\exp(c_2 x))\mathrm{d}x
 \\ 
 = \: & 
 \frac{2C_0}{c_2}\exp(-c_1\exp(c_2 Nh)).
\end{align*}
Then, we obtain (\ref{eq:error-estimate}).
\hfill\rule{1.5ex}{1.5ex}

\medskip

We remark here that we can reduce the number of sampling points by half if the integrand $f(x)$ is real valued 
on the real axis. 
In fact, we have 
$f(\overline{z}) = \overline{f(z)}$ 
by the reflection principle, taking the integral path $C$ symmetric with respect to the real axis, 
that is, 
$\varphi(-u) = \overline{\varphi(u)}$, 
which leads to 
\begin{equation*}
 \varphi^{\prime}(-u) = - \overline{\varphi^{\prime}(u)}, 
\end{equation*}
and taking the DE transform $\psi_{\rm DE}(v)$ to be an even function, 
we have
\begin{align}
% \nonumber
% & 
% \fp\int_0^{\infty}x^{-n}f(x)\mathrm{d}x
% \\
 \nonumber
 & 
 I^{(n)}[f]
 \simeq 
 I_{h,N}^{(n)\prime}[f]
 \\
 \nonumber
 \equiv \: & 
 \frac{h}{2\pi}\im
 \left\{
 \varphi(\psi_{\rm DE}(0))^{-n}f(\varphi(\psi_{\rm DE}(0)))
 \log(-\varphi(\psi_{\rm DE}(0)))\varphi^{\prime}(\psi_{\rm DE}(0))\psi_{\rm DE}^{\prime}(0)
 \right\}
 \\
 \nonumber
 & 
 + 
 \frac{h}{\pi}\im
 \bigg\{
 \sum_{k=1}^{N}\varphi(\psi_{\rm DE}(kh))^{-n}
 f(\varphi(\psi_{\rm DE}(kh)))\log(-\varphi(\psi_{\rm DE}(kh)))
 \\ 
 & 
 \label{eq:approx-fp-integral2}
 \hspace{15mm}
 \times 
 \varphi^{\prime}(\psi_{\rm DE}(kh))\psi_{\rm DE}^{\prime}(kh)
 \bigg\}.
\end{align}
%%%%%%%%%%%%%%%%%%%%%%%%%%%%%%%%%%%%%%%%%%%%%%%%%%%%%%%%%%%%%%%%%%%%%%%%%%%%%%%%%%%%%%%%%%%%%%%%%%%%%
\section{Numerical examples}
\label{sec:example}
In this section, we show some numerical examples which show the effectiveness of the proposed method.

\begin{figure}[htbp]
 \begin{center}
  \psfrag{x}{$\re z$}
  \psfrag{y}{$\im z$}
  \includegraphics[width=0.7\textwidth]{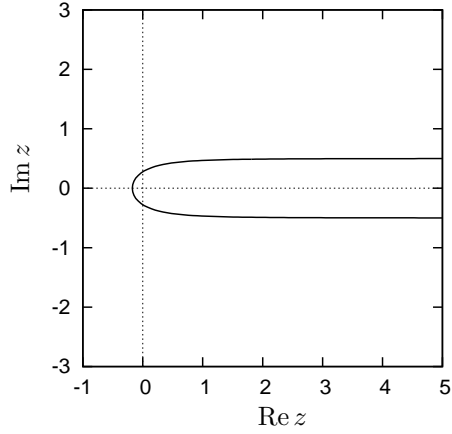}
 \end{center}
 \caption{The complex integral path $C$.}
 \label{fig:example-integral-path}
\end{figure}

\begin{figure}[htbp]
 \begin{center}
  \begin{tabular}{cc}
   \psfrag{N}{$N$} 
   \psfrag{e}{\rotatebox{-90}{\hspace{-8mm}$\varepsilon_{N}^{(n)}[f]$}}
   \includegraphics[width=0.45\textwidth]{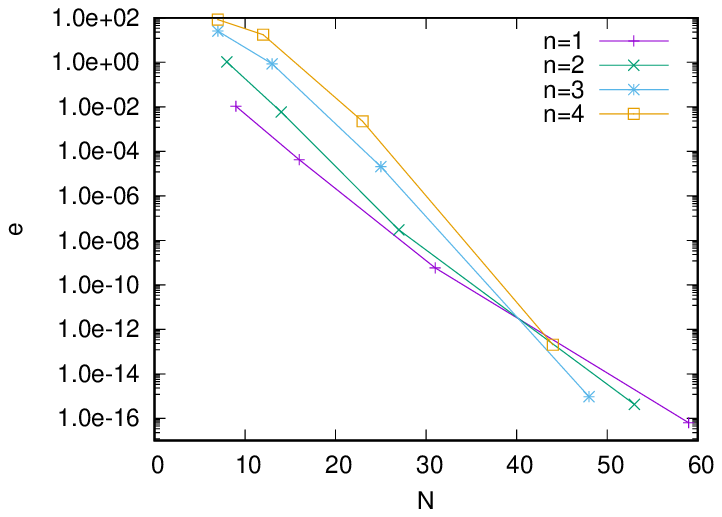}
   & 
       \psfrag{N}{$N$}
       \psfrag{e}{\rotatebox{-90}{\hspace{-8mm}$\varepsilon_{N}^{(n)}[f]$}}
       \includegraphics[width=0.45\textwidth]{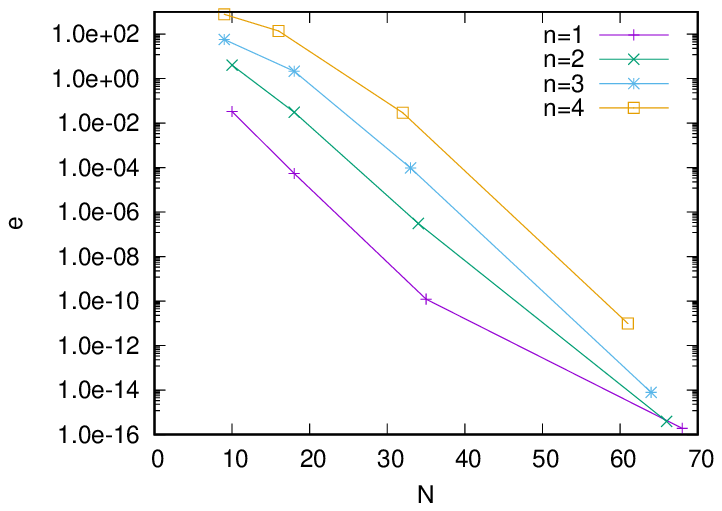}
       \\
   integral (i) & integral (ii)
  \end{tabular}
 \end{center}
 \caption{The errors of the proposed approximation formula (\ref{eq:approx-fp-integral2}) 
 applied to the f.p. integrals (\ref{eq:example}).}
 \label{fig:example}
\end{figure}

We computed the f.p. integrals
\begin{equation}
 \label{eq:example}
  \begin{aligned}
   \mathrm{(i)} \quad & \fp\int_0^{\infty}\frac{x^{-n}}{1 + x^2}\mathrm{d}x 
   = 
   \begin{cases}
    (\pi/2)(-1)^m & ( \: n = 2m \ \mbox{(even)} \: ) \\ 
    0 & ( \: n = 2m+1 \ \mbox{(odd)} \: )
   \end{cases}
   \\ 
   \mathrm{(ii)} \quad & \fp\int_0^{\infty}x^{-n}\mathrm{e}^{-x}\mathrm{d}x
   = 
   \begin{cases}
    - \gamma & ( \: n = 1 \: ) \\ 
    - 1 + \gamma & ( \: n = 2 \: ) \\ 
    \frac{3}{4} - \frac{1}{2}\gamma & ( \: n = 3 \: ) \\ 
    - \frac{11}{36} + \frac{1}{6}\gamma & ( \: n = 4 \: )
   \end{cases}
  \end{aligned}
\end{equation}
for $n=1,2,3,4$, where $\gamma$ is Euler's constant,  
by the formula (\ref{eq:approx-fp-integral2}). 
All the computations were performed using programs coded in C++ with double precision working. 
The complex integral path $C$ in (\ref{eq:complex-integral}) is taken as 
\begin{equation*}
 C : \: z = \varphi(u) = 
  \frac{u+0.5\mathrm{i}}{\mathrm{i}\pi}\log\left(\frac{1+\mathrm{i}(u+0.5\mathrm{i})}{1-\mathrm{i}(u+0.5\mathrm{i})}\right), 
  \quad +\infty > u > -\infty 
\end{equation*}
(see Figure \ref{fig:example-integral-path}). 
We took the number of sampling points $N$ for given mesh $h=2^{-1}, 2^{-2}, \ldots $
by truncating the infinite sum at the $k$-th term such that 
\begin{equation*}
 \frac{h}{\pi}\times|\mbox{the $k$-th term}| < 
  \begin{cases}
   10^{-15}\times |I_{h,N}^{(n)\prime}[f]| & \mbox{if} \quad I_{h,N}^{(n)\prime}[f] \neq 0 
   \\ 
   10^{-15} & \mbox{otherwise}.
  \end{cases}
\end{equation*}
Figure \ref{fig:example} shows the relative errors of the proposed approximation formula (\ref{eq:approx-fp-integral2}) 
\begin{equation*}
 \varepsilon_{N}^{(n)}[f] = 
  \begin{cases}
   |I_{h,N}^{(n)}[f] - I^{(n)}[f]| / |I^{(n)}[f]| & ( \: I^{(n)}[f] \neq 0 \: ) \\ 
   |I_{h,N}^{(n)}[f]| & \mbox{(otherwise)}
  \end{cases}
\end{equation*}
applied to the f.p. integrals (\ref{eq:example}). 
These figures shows the exponential convergence of the proposed formula as the number of sampling points $N$ increases. 
%%%%%%%%%%%%%%%%%%%%%%%%%%%%%%%%%%%%%%%%%%%%%%%%%%%%%%%%%%%%%%%%%%%%%%%%%%%%%%%%%%%%%%%%%%%%%%%%%%%%%
\section{Summary}
\label{sec:summary}
In this paper, we proposed a numerical integration formula for Hadamard finite-part integrals 
with an integral power singularity at the endpoint on a half-infinite interval. 
In the proposed method, we express the desired f.p. integral using a complex integral, 
and we obtain the f.p. integral by evaluating the complex integral by the DE formula. 
Theoretical error estimate and some numerical examples show the exponential convergence of the proposed method 
in the case that the integrand is an analytic function. 
%%%%%%%%%%%%%%%%%%%%%%%%%%%%%%%%%%%%%%%%%%%%%%%%%%%%%%%%%%%%%%%%%%%%%%%%%%%%%%%%%%%%%%%%%%%%%%%%%%%%%
\bibliographystyle{plain}
\bibliography{arxiv2019_4}

\begin{thebibliography}{10}

\bibitem{Bialecki1990a}
B.~Bialecki.
\newblock A sinc-hunter quadrature rule for cauchy principal value integrals.
\newblock {\em Math. Comput.}, 55:665--681, 1990.

\bibitem{Bialecki1990b}
B.~Bialecki.
\newblock A sinc quadrature rule for hadamard finite-part integrals.
\newblock {\em Numer. Math.}, 57:263--269, 1990.

\bibitem{ElliotPaget1979}
D.~Elliot and D.~F. Paget.
\newblock Gauss type quadrature rules for cauchy principal value integrals.
\newblock {\em Math. Comput.}, 33:301--309, 1979.

\bibitem{EstradaKanwal1989}
R.~Estrada and R.~P. Kanwal.
\newblock Regularization, pseudofunction, and hadamard finite part.
\newblock {\em J. Math. Anal. Appl.}, 141:195--207, 1989.

\bibitem{Ogata2019c}
H.~Ogata.
\newblock A numerical method for computing hadamard finite-part integrals with
  a non-integral power singularity at an endpoint, 2019.
\newblock arXiv:1909.11398v1 [math.NA].

\bibitem{Ogata2019b}
H.~Ogata.
\newblock A numerical method for hadamard finite-part integrals with an
  integral power singularity at an endpoint, 2019.
\newblock arXiv:1909.08872v1 [math.NA].

\bibitem{OgataHirayama2018}
H.~Ogata and H.~Hirayama.
\newblock Numerical integration based on hyperfunction theory.
\newblock {\em J. Comput. Appl. Math.}, 327:243--259, 2018.

\bibitem{OgataSugiharaMori2000}
H.~Ogata, M.~Sugihara, and M.~Mori.
\newblock De-type quadrature formulae for cauchy principal-value integrals and
  for hadamard finite-part itnegrals.
\newblock In {\em Proceedings of the Second ISAAC Congress}, volume~1, pages
  357--366, 2000.

\bibitem{Paget1981}
D.~F. Paget.
\newblock The numerical evaluation of hadamard finite-part integrals.
\newblock {\em Numer. Math.}, 36:447--453, 1981.

\bibitem{Stenger1993}
F.~Stenger.
\newblock {\em Numerical Methods Based on Sinc and Analytic Functions}.
\newblock Springer-Verlag, New York, 1993.

\bibitem{TakahasiMori1974}
H.~Takahasi and M.~Mori.
\newblock Double exponential formulas for numerical integration.
\newblock {\em Publ. Res. Inst. Math. Sci., Kyoto Univ.}, 339:721--741, 1978.

\end{thebibliography}
\end{document}